\input amstex
\documentstyle{amsppt}
\magnification=1200 
\hcorrection{.25in}
\topmatter
\title  The Krein spectral shift and \hbox{rank one} perturbations of spectra 
\endtitle
\author A. Poltoratski\thanks Research at MSRI is supported in part
by NSF grant DMS-9022140 \endthanks\endauthor
\address MSRI, 1000 Centennial Drive, Berkeley, California 94720\endaddress
\email agp\@msri.org \endemail
\abstract We use recent results on the boundary behavior of 
Cauchy integrals to study the
Krein spectral shift of a rank one perturbation problem for
 self-adjoint operators. As an application, we prove that
all self-adjoint rank one 
perturbations of a self-adjoint operator are pure point
if and only if the spectrum of the operator is countable.
We also study pairs of pure point operators unitarily equivalent up to a 
rank one perturbation and give various examples of rank one
perturbations of singular spectra.
\endabstract
\endtopmatter
\document

\tolerance=1000

\mathchardef\lesha="7F1F

\pageno1

\define\nec{\underset{ z @>>\ngtr> x }\to                       
{-\kern -5pt -\kern -5pt-\kern -5pt-\kern -5pt\longrightarrow}}
\define\Clos{\operatorname{Clos}}
\define\supp{\operatorname{supp}}
\define\im{\operatorname{Im}} 
\define\re{\operatorname{Re}}
\define\dist{\operatorname{dist}}

\define\a{\Cal A}
\define\Cb{\Cal B}\define\Cc{\Cal C}

\heading 1. Introduction. \endheading

In this paper 
we  use   well-known connections
between one-dimensional perturbation theory
and complex analysis. We apply recent results
on the boundary behavior of Cauchy integrals
to study the spectral
properties of rank one perturbations of self-adjoint operators. 

 Section 2  discusses the structure of a family of the spectral measures
of rank one perturbations of a self-adjoint operator. 
We establish the relations between such families and the 
characteristic function of the operator. 

  Section 3 is devoted to
the properties of the Krein spectral shift of a rank one perturbation 
problem. We study spectral shifts 
in terms of the boundary behavior of the resolvent function. 

In Section 4 we prove that 
all self-adjoint rank one 
perturbations of a self-adjoint operator are pure point
if and only if the spectrum of the operator is countable.

 Section 5 analyzes relative spectral properties
of two self-adjoint pure point operators
unitarily equivalent up to a rank one perturbation.

In Section 6 we give an example of 
``the absence of mixed spectra''. We construct
a family $A+\lambda(\cdot,\varphi)\varphi$  of rank one perturbations 
of a self-adjoint operator such that the operators corresponding
to the coupling constants $\lambda\in [0;1] $ are singular continuous
and the operators corresponding to the coupling constants
$\lambda\in \Bbb R\setminus [0;1]$
are pure point.

\remark{\bf Acknowledgments} I would like to thank Nikolai Makarov for asking some of the questions
discussed in this paper. I am also grateful to Aleksei Aleksandrov, Dmitri Jakubovich and
Fedor Nazarov for useful discussions and  comments, and to the Mathematical
Sciences Research Institute for hospitality while this work was completed.
\endremark

 \heading  2. Families $\Cal M_\varphi$ and rank one perturbations  
 of self-adjoint operators.
 \endheading               
 We will denote by  $\Cal M( \hat\Bbb R)$ the space of Borel 
 complex measures on $ \hat\Bbb R=\Bbb R\cup{\infty}$ with the norm
 $$||\mu||=\frac 1\pi\int\limits_{-\infty}^\infty \frac {d|\mu(t)|}{1+t^2}+|\mu|(\infty).$$
 We will also use the notation $\Cal M_+( \hat\Bbb R)$ for the subset
 consisting of nonnegative measures.

  Let $\phi$  be a nonconstant analytic function  in the upper half plane $\Bbb C_+$ 
  such that $|\phi|\leq 1$.
  Then  for any $\alpha$ from the unit circle $\Bbb T$  
 function $\frac{\alpha+\phi}{\alpha-\phi}$ has positive real part
in $\Bbb C_+$ . Thus there exists
  a measure $\mu_\alpha\in \Cal M_+( \hat\Bbb R)$ 
  such that its Poisson integral satisfies
$$ (\Cal P\mu_\alpha)(x+iy)=
 \frac 1\pi \int\limits_{-\infty}^{\infty}
\frac{y d\mu_\alpha(t)}{(x-t)^2+y^2} +y\mu_\alpha(\infty)
=\re
\frac{\alpha+\phi}{\alpha-\phi}
 .\tag 2.1$$
We will denote by $\Cal M_\phi( \hat\Bbb R)$
the family of such measures $\{\mu_\alpha\}_{\alpha\in\Bbb T}\ $ 
corresponding to function $\phi$. 

If $\phi$ is defined in the unit 
disk one can replace the Poisson integral in (2.1) by the Poisson integral in the disk
and consider the analogous family $\Cal M_\phi(\Bbb T)$ consisting of
positive measures on the unit circle $\Bbb T$. 

Families  
$\Cal M_\phi $ 
possess many interesting properties, see [A1], [A2], [C], [P1], [P2] and  [Sa].
 As  was shown by Clark in [C], 
  if $\phi$ is an inner
function in the unit disk $\Bbb D$  then 
$\Cal M_\phi(\Bbb T)$ is the system of the spectral measures
of all unitary rank one  perturbations of 
the model contraction with the characteristic function $\phi$.
By adopting
the argument of Clark to the case of the upper half-plane, one can make a 
similar connection  between  
 rank one perturbations of self-adjoint operators and families 
 $\Cal M_\phi( \hat\Bbb R)$. Here we  present a different way to establish
this connection.

If  $A_0$ is a bounded cyclic self-adjoint operator
and  $\varphi$ is its cyclic vector, then we can consider the family of
 rank one  perturbations $$A_\lambda=A_0+
\lambda(\cdot , \varphi)\varphi , $$ $\lambda\in\Bbb R  . $
Let  $\nu_\lambda$ 
be the spectral measure of $\varphi$ for $A_\lambda$. Since
$A_\lambda$ is bounded, $\supp\nu_\lambda\subset\Bbb R$ is compact. Thus 
the Cauchy transform  of
$\nu_\lambda$

 $$\Cal K\nu_\lambda(z)=\frac 1\pi\int\limits_{-\infty}^{\infty}\frac {d\nu_
 \lambda(t)}{t-z}=\frac 1\pi((A_\lambda-z )^{-1}\varphi , \varphi)\tag 2.2$$
converges for any $z\in\Bbb C_+$.

The relation for the resolvents$$
 (A_0 - z)^{-1} -(A_\lambda-z)^{-1}
 =\left[\lambda(\cdot , \varphi)\left((A_\lambda
 -z)^{-1}\varphi\right)\right]
 (A_0-z)^{-1},$$
for any $z\in\Bbb C_+$, 
gives us $$
  \Cal K\nu_\lambda(z)=\frac {\Cal K\nu_0(z)}{1+\pi\lambda \Cal K\nu_0(z)},\tag 2.3 $$ 
see [A].

Since $\nu_0$ is a nonnegative
Borel measure  on $ {\Bbb R}$, 
$-i\Cal K\nu_0$ is  an analytic function with a nonnegative  real part in $\Bbb C_+$:
$$\re -i\Cal K\nu_0=\Cal P\nu_0.$$ 
Thus $$-i\Cal K\nu_0=\frac{1+\phi}{1-\phi}\tag 2.4$$ 
for some bounded analytic 
function $\phi,\  |\phi|\leq1$ in $\Bbb C_+$.
 
Hence by (2.3) 
 $$ (\Cal P\nu_\lambda)(x+iy)=\re -i\Cal K\nu_\lambda(x+iy)=\re \frac{\frac{1-\phi}{1+\phi}}
{1+i\pi\lambda\frac{1-\phi}{1+\phi}}
= c\re\frac{\beta-\phi}{\beta+\phi}
\tag2.5$$
 where     $\beta=\frac{1+i\pi\lambda}{1-i\pi\lambda} $,
$c=\frac1{1+\pi^2\lambda^2} $. 
Thus
$$ \Cal P\nu_\lambda=c\Cal P\mu_\beta$$ 
where $\{\mu_\beta\}_{\beta\in\Bbb T}=\Cal M_\phi( \hat\Bbb R)$.

Conversely, if
$\phi$ is an analytic function 
in $\Bbb C_+$ with $|\phi|\leq 1$ we can consider
the operator $A_\mu$ of multiplication by $ z$ in $L^2(\mu)$
where $\mu=\mu_1\in\Cal M_\phi( \hat\Bbb R)$. 
Reversing the above argument
we can show that for any $\alpha\in\Bbb T$ the measure $\mu_\alpha\in\Cal M_\phi$ is
a spectral measure  for the rank one perturbation $A+\lambda(\cdot,1)1$ where
$\lambda=\frac i\pi\frac {1-\alpha}{1+\alpha}$.
Therefore, each family $\Cal M_\phi( \hat\Bbb R)$ is a family of spectral measures of 
rank one perturbations of some self-adjoint operator.

Each measure from $\Cal M_+(\hat\Bbb R)$ belongs to a family $\Cal M_\phi$ for some $\phi$
($\phi$ can be found from (2.1) uniquely up to a M\"obius transform).
One can also say that two given measures $\mu$ and $\nu$  from $\Cal M(\hat\Bbb R)$ belong
to the same family $\Cal M_\phi$ iff there exist operator $A$ and its rank one perturbation
$A_\lambda=A+\lambda(\cdot,\varphi)\varphi$ such that $\mu$ and $\nu$ are the spectral measures of $\varphi$
for $A$ and $A_\lambda$. Indeed,
as we proved above, such $A$ and $A_\lambda$ exist iff 
$c_1\mu,c_2\nu\in\Cal M_\psi$
for some function $\psi$ and positive constants $c_1$ and $c_2$. 
Suppose $c_1\mu=\mu_\alpha, c_2\nu=\mu_\beta$ for some $\mu_\alpha,\mu_\beta\in\Cal M_\psi$.
Then
$\mu,\nu\in\Cal M_\phi$ for $\phi=b\circ\psi$ where $b$ is the M\"obius transform of the unit disk
 $\Bbb D$
 such that $|b'(\alpha)|=1/c_1, |b'(\beta)|=1/c_2$, see [A2].

Further relations between the resolvent functions of the rank one perturbations
of a self-adjoint operator are  discussed in
 [A],  [D], [M-P],  [S-W], [S] and [R-J-L-S].

\remark{\bf Remark}
The same argument works for the families of unitary operators
(in particular for those studied by Clark [C]).

Let $\Cal U_1$ be a unitary cyclic operator. 
Let probability measure $\mu_1\in\Cal M_+(\Bbb T)$ be 
 the spectral measure of some cyclic vector $v$ for $\Cal U_1$. 
 Then we can consider the family of
unitary rank one perturbations of $\Cal U_1$ : $$\Cal U_\alpha=
\Cal U_1+(\alpha-1)(\cdot , \Cal U_1^{-1}v)v, $$ $\alpha\in\Bbb T$. For the resolvents we have 
$$(\Cal U_1-z)^{-1}-(\Cal U_\alpha-z)^{-1}=
(\Cal U_1-z)^{-1}[\alpha(\cdot , \Cal U_1^{-1}v)v](\Cal U_\alpha-z)^{-1}
.$$ 

Denote by $\mu_\alpha$ the spectral measure of $v/||v||$ for
$U_\alpha$. Let $ K\mu_\alpha$ and $ P\mu_\alpha$ be the standard Cauchy and Poisson 
integrals of measure $\mu_\alpha$ in the unit disk $\Bbb D$. Since
$ K\mu_\alpha=((\Cal U_\alpha-z)^{-1}v , v)$ for $z\in\Bbb D$,
we have that
$$ K\mu_\alpha=\frac{\alpha K\mu_\alpha}
              {1+ (\alpha-1) K\mu_1} .\tag 2.6 $$

If we consider an analytic function
 $\theta$  in $\Bbb D$  such that $ |\theta|\leq1, \theta(0)=0$
 and $$ P\mu_1=\re\frac{1+\theta}{1-\theta}$$ 
in $\Bbb D$, then 
$$ K\mu_1=\frac1{1-\theta}$$
and by (2.6)
$$ P\mu_\alpha=2\re  K\mu_\alpha-1=
2\re\left[\frac{\frac\alpha{1-\theta}}{1+\frac{\alpha-1}{1-\theta}}-\frac12
\right]=
\re \frac{\alpha+\theta}{\alpha-\theta}.$$
 
 Thus $\{\mu_\alpha\}_{\alpha\in\Bbb T}
=\Cal M_\theta (\Bbb T)$. 
\endremark

As we mentioned before,
this result was obtained  in
[C] for 
 rank one  unitary
perturbations of the model contraction 
$T_\theta=\Cal S P_\theta $ 
 , where 
 $\theta$ is an inner function, 
$\Cal S:f\mapsto zf$ is a shift operator 
in the Hardy space $ H^2$  and $ P_\theta$
is the orthogonal projector from $ H^2$ onto the model space  
$\theta^*( H^2)=
 H^2\circleddash\theta H^2$. Any singular cyclic unitary operator
can be represented as a rank one perturbation of $T_\theta$ for
some inner $\theta$. Conversely, if the operator $\Cal U_1$ from the last remark
is singular, then $\  \Cal U_1-(\cdot , \Cal U_1^{-1}v)v $ 
is a $C_0$ completely nonunitary
contraction with the characteristic function $\theta$.

 \heading 3. The Krein spectral shift.   \endheading

In this section we discuss the notion of the Krein spectral shift
in the case of rank one perturbations.
For basic results and definitions in this area we refer to [K], [A-D], [M-P] and [S].

Let $A$ be a bounded self-adjoint operator, $\varphi $ its cyclic vector.
Let $\mu\in\Cal M_+(\hat\Bbb R)$ be the spectral  measure of $\varphi$ for $A$.
Since $\mu$ is a nonnegative measure, its Cauchy transform has nonnegative
imaginary part. Therefore for each $\lambda\in\Bbb R$ the function $1+\pi\lambda\Cal K$
admits the following representation:
$$1+\pi\lambda\Cal K\mu=\exp[\Cal Ku]\tag 3.1$$ 
fore some real valued function $u\in L^\infty(\Bbb R), || u||\leq\pi$ with compact support.
If $\lambda\geq 0$ then function $u$ can be chosen to satisfy 
$$u=\arg(1+\pi\lambda\Cal K\mu)\tag 3.2$$
where $\arg$ stands for the principal branch of argument taking values
in $(-\pi;\pi]$
(we assume that $\Cal K\mu$ is defined a. e.
on $\Bbb R$ by its angular boundary values). In this case we will have
$0\leq u\leq \pi$. Similarly
if $\lambda<0$  then $u$ can be chosen to satisfy (3.2) with
$\arg$ taking values in $[-\pi;\pi)$. In this case we will have $0\leq u\leq -\pi$. 
\definition{Definition}
Such $u$  is called the 
Krein spectral shift of the
perturbation problem $\left(A\mapsto A+\lambda(\cdot,\varphi)\varphi\right)$.

\enddefinition

The Krein spectral shift is usually defined in more general settings of
compact perturbations $\left(A\mapsto A+K\right)$, see [K] and [M-P].
 When $K$ is one-dimensional
 the general definition can be reduced to the one given above, see [M-P].

 Let $\nu$ be the spectral measure of $ \varphi$ for $A+\lambda(\cdot,\varphi)\varphi$.
Even though measure $\nu$ does not appear in (3.1), it is uniquely determined by $\mu$ and $\lambda$ and 
satisfies a similar equation.

\proclaim{Lemma 3.1}

Measures $\mu$, $\nu$ and the Krein spectral shift $u$ satisfy
$$1+\pi\lambda\Cal K\mu=\exp[\Cal K(u)]=\left[1-\pi\lambda\Cal K\nu\right]^{-1}.\tag3.3$$ 
\endproclaim
\demo{Proof} Apply formulas (3.1) and 
 (2.3) with $\nu_0=\mu$ and $\nu_\lambda=\nu$. 
$\blacktriangle$
\enddemo

Since in this paper we will mostly deal with the spectral measures  of operators,
it will be more convenient for us to call $u$
the phase shift of the pair of measures $(\mu,\nu)$ for the coupling constant $\lambda$.
If $\lambda=1$ we will call $u$ the phase shift of the pair $(\mu;\nu)$. One can use formula (3.3) 
 as an alternative definition of the phase shift.

Note that in this definition $u$ determines $\mu$ and $\nu$ uniquely. 
This reflects the fact that among the  spectral measures of $A$ and 
$A+\lambda(\cdot,\varphi)\varphi$ we always choose those corresponding
to the vector $\varphi$.

Any function 
$u\in L^\infty(\Bbb R), 0\leq u \leq\pi$ with compact support
is the Krein spectral shift of some perturbation problem (the phase shift of some pair of measures). Indeed,
for  any $\lambda>0$  one can consider a nonnegative measure
$\mu\in\Cal M(\hat\Bbb R)$ with compact support satisfying (3.1).
Then $u$ 
is the Krein spectral shift of the perturbation problem 
$\left(A_\mu\mapsto A_\mu+\lambda(\cdot,1)1\right)$. 
Similarly each
$u\in L^\infty(\Bbb R), 0\geq u \geq-\pi$ with compact support
is also the Krein spectral shift of some  perturbation problem (with $\lambda<0$).

As was shown in the  previous Section, there exist
operator $A$ and its cyclic vector $\varphi$ 
such that $\mu$ and $\nu$ are spectral measures
for $\varphi$ of $A$ and $A+(\cdot,\varphi)\varphi$ 
iff both measures belong to the same
 family $\Cal M_\phi$ for some function $\phi$ analytic in $\Bbb C_+$.
One can show that if  we take $\mu_\alpha,\mu_\beta\in\Cal M_\phi$ 
$\alpha=e^{i\theta_1} ,\beta=e^{i\theta_2}, 0\leq \theta_1<\theta_2<2\pi$ where $\phi$
is inner then the phase shift $u$ for the pair   
$(\mu_\alpha,\mu_\beta)$ can be defined as  $$u=\pi\chi_{\{x\in\Bbb R|\phi(x)=e^{i\theta}, \theta_1<\theta<\theta_2\}}$$
where $\chi_E$ denotes the indicator function of the set $E$ (we again assume that $\phi$ is defined a. e.
on $\Bbb T$ by its angular boundary values).

As before, let $u$ be the phase shift of the pair of measures $(\mu,\nu)$ and $\phi$ be a characteristic function:
$\mu,\nu\in\Cal M_\phi$. 
Then $u$ has the following properties.

\proclaim{Lemma 3.2}

Let $K\subset \Bbb R$ be a measurable set.
 Then the following conditions are equivalent 

i) the restrictions of $\mu$ and $\nu$ on $K$ are singular,

ii) $|u|$ is equal to $0$ or $\pi$ a. e. on $K$,

iii) $|\phi|=1$ a. e. on $K$.

\endproclaim
\demo{Proof} For i)$\Leftrightarrow $ii) see [M-P] or [S];
for i)$\Leftrightarrow$iii) see [C].$\blacktriangle$
\enddemo

We will denote by $Qu$ the conjugate Poisson integral of $u$:
$$Qu(x+iy)=-\re \Cal Ku(x+iy)=\frac1\pi\int\limits_{-\infty}^{\infty}
\frac{x-t}{(x-t)^2+y^2}u(t)dt.$$
It is well-known that if $v\in L^\infty (\Bbb R)$ with a compact support then the difference 
$$ Qv(x+iy)- \int\limits_{\Bbb R\setminus (x-y;x+y)}\frac {v(t)dt}
{x-t}\tag 3.4$$
is $O(1)$
as $x+iy\underset\ngtr\to\rightarrow x_0$ for any $x_0\in\Bbb R$ and $o(1)$ for any Lebesgue point $x_0$
of $v$, see [G].
 
Let $E\subset\Bbb R$ be a measurable set, $\mu\in \Cal M(\Bbb R)$.
We will denote
$$p.\ v.\ \int\limits_Ed\mu(x)=
\lim_{\epsilon\to 0}\int\limits_{E\setminus
(-\epsilon;\epsilon)}d\mu(x).$$

We will denote by $\mu^s$ the singular component of measure $\mu$.
The definition of the phase shift and relation (3.4) imply the following

\proclaim{Lemma 3.3}

i)
 For $\mu^s$-a.e $x$ 
$$p.\ v.\ \ \int\limits_\Bbb R \frac{u(x+t)dt}{t}= \infty.\tag 3.5$$
 For $\nu^s$-a.e $x$ 
$$p.\ v.\ \ \int\limits_\Bbb R \frac{u(x+t)dx}{t}=- \infty.\tag 3.6$$

ii)    
There exist $\alpha,\beta\in\Bbb T$ such that
 for $\mu^s$-a.e $x$ $$\phi(z) \nec \alpha;$$
 For $\nu^s$-a.e $x$ $$\phi(z) \nec \beta.$$
\endproclaim
\demo{Proof}Follows from  formulas (2.1) and (3.3). 
For i) see [S]; for ii) see [C].$\blacktriangle$
\enddemo

\example{\bf Example 3.4}
Let $x,y\in\Bbb R, x<y$. Let
 $u=\pi \chi_{(x;y)}$. Then by Lemma 3.2 ii) $\mu$ and $\nu$ are singular. By Lemma 3.3 i)
$\mu$ is concentrated on the set $\{x\}$ and $\nu$ is concentrated on the set 
$\{y\}$. Since $u$ is nonconstant, $\mu$ and $\nu$ are nonzero.
Thus $\mu=b\delta_{x}$ and  $\nu=c\delta_{y}$ for some positive $b$ and $c$.
\endexample

To prove the next property of phase shifts we will need the following

\proclaim{Theorem 3.5([P1])}
Let $\sigma\in M( \hat\Bbb R)$ be a singular measure with compact support,
$f\in L^1(|\sigma|)$. 
Define  meromorphic function $F$ in $\Bbb C_+$ as
$$F=\frac{1+\Cal K(f\sigma)}{1+\Cal K\sigma}.$$
Then $F$ has  nontangential boundary limits equal to $f$$\ \ $ 
$|\sigma|$-a. e. 
on $ {\Bbb R}$.
\endproclaim

This Theorem is proved in [P1] in the case of the unit disk.

We will say that real function $f$ defined in $\Bbb C_+$ is less than $\infty$
(greater than $-\infty$) at point $x\in\Bbb R$ if $$\liminf_{z\underset \ngtr \to \longrightarrow x}
f(z)<\infty \ \ \  (\limsup>-\infty).$$
We will write that $$p.\ v.\ \int\limits_E f(x)dx<\infty \ \ (>-\infty)$$
if $$\liminf_{\epsilon\to 0}\int\limits_{E\setminus
(-\epsilon;\epsilon)}d\mu(x)<\infty\ \ \ (\limsup>-\infty).$$
	
\proclaim{Corollary 3.6} Let $\sigma,\gamma\in\Cal M_+{\Bbb R})$ be  singular measures with compact supports.
Define  holomorphic function $F$ in $\Bbb C_+$ as
$$F=\frac{1+\Cal K\sigma}{1+\Cal K\gamma}.$$
Let $K$ be a measurable subset of $\Bbb R$. Then $|F|<\infty$ $\sigma$-a. e. on $K$
iff the restriction of $\sigma$ on $K$ is absolutely continuous with
respect to $\gamma$.
\endproclaim
\demo{Proof}
Let $\sigma=f\gamma+\eta$ where $f\geq 0, f\in L^1(\gamma),\eta\perp\gamma$.
Then
 $$F=\frac{
1+\Cal K\sigma}{1+\Cal K\gamma}=\frac{1+\Cal K\sigma}{1+\Cal K(\sigma+\gamma)}\times
\frac{1+\Cal K(\sigma+\gamma)}{1+\Cal K\gamma}.\tag 3.7$$
By Theorem 3.5 the first fraction has finite nonzero limits  $|\sigma|$-a. e.
Since
$$
\frac{1+\Cal K(\sigma+\gamma)}{1+\Cal K\gamma}=\left[\frac{1+\Cal K\gamma}{1+\Cal K((1+f)\gamma+\eta)}\right]^{-1},$$
by Theorem 3.5  the second fraction in (3.7) has finite limits $|\gamma|$-a. e. and 
tends to infinity $\eta$-a. e. Thus  $|F|<\infty$ $\sigma$-a. e.
on $K$ iff $\eta(K)=0$.
$\blacktriangle$
\enddemo

\proclaim{Lemma 3.7}
Let $u_1$ and $u_2$ be the  phase shifts of the pairs of measures
$(\mu_1,\nu_1)$ and $(\mu_2,\nu_2)$ respectively. Let $K\subset\Bbb R$
be a measurable set.
Then 

i)  
$$p.\ v.\ \int\limits_\Bbb R (u_1(t+x)-u_2(t+x))\frac{dt}{t}<\infty$$
for  $\mu_1$-a. e.  $x\in K$ iff the restriction
of $\mu_1^s$ on $K$ is absolutely continuous with respect to
$\mu_2$; 
$$p.\ v.\ \int\limits_\Bbb R (u_1(t+x)-u_2(t+x))\frac{dt}{t}>-\infty$$
for  $\nu_1$-a. e.  $x\in K$ iff 
the restriction
of $\nu_1^s$ on $K$ is absolutely continuous with respect to
$\nu_2$;

ii) if 
$\mu_2^s$-a. e. $x\in K$ is a Lebesque point of $u_1-u_2$ and 
$$p.\ v.\ \int\limits_\Bbb R (u_1(t+x)-u_2(t+x))\frac{dt}{t}= f(x)<\infty$$
 for 
$\mu_2^s$-a. e. $x\in K$ then
the restriction
of $\mu_1^s$ on $K$ is equal to the restriction of $e^{f}\mu_2^s$ on $K$;
 if 
$\nu_2^s$-a. e. $x\in K$ is a Lebesque point of $u_1-u_2$ and 
$$p.\ v.\ \int\limits_\Bbb R (u_1(t+x)-u_2(t+x))\frac{dt}{t}= f(x)>-\infty$$
for
$\nu_2^s$-a. e. $x\in K$ then
the restriction
of $\nu_1^s$ on $K$ is equal to the restriction of $e^{-f}\nu_2^s$ on $K$.

\endproclaim
\demo{Proof}
Notice that by the defining formula (3.3)
$$\exp\left[\Cal K(u_1-u_2)\right]=\frac{1+\Cal K\mu_1}{1+\Cal K\mu_2}
=\frac{1+\Cal K\nu_2}{1+\Cal K\nu_1}.$$
Now   Corollary 3.6 and formula (3.4) imply part $i)$. Since the difference (3.4) is $o(1)$ at Lebesgue
points of $u_1-u_2$, part $ii)$
follows from part $i)$ and Theorem 3.5.$\blacktriangle$
\enddemo

\proclaim{Lemma 3.8 [M-P]} Measure  $\mu$ has a point mass  at                    
$x$ iff
$$\int\limits_{x-1}^{x+1} (\pi\chi_{(x;x+1)}-u(y))\frac {dy}{y-x}<\infty.\tag3.8$$
Measure
 $\nu$ has a point mass  at                    
$x$ iff
$$\int\limits_{x-1}^{x+1} (\pi\chi_{(x-1;x)}-u(y))\frac {dy}{x-y}<\infty.\tag3.9$$
\endproclaim
\demo{Proof}
As was shown in Example 3.5 measure $\mu_0$ corresponding to the phase shift
$u_1=\pi\chi_{(x;x+1)}$ has a point mass at $x$. It is left to apply statement i) of  Lemma 3.7
to $u_1, u_2=u$, and $K=\{x\}$. 

The statement can also be proved directly, see [M-P].
$\blacktriangle$\enddemo

\remark{\bf Remark} In terms of the characteristic function $\phi$ a point mass
of $\mu_\alpha\in\Cal M_\phi$ can be described in the following way, see [C], [A-C] or [Sa].
For any $x\in\Bbb R$ $\mu_\alpha(x)>0$ iff 
$\phi(z)\to \alpha
$ as $z \underset \ngtr \to \longrightarrow x$ and
$\phi$ has a finite nontangential derivative $\phi'(x)$ at $x$.

In terms of other measures from the family $\Cal M_\phi$, $\mu_1(x)>0$ iff
$\Cal K\mu_{-1}(z)\to 0$ as $z\underset \ngtr \to \longrightarrow x$ and
$(x-y)^{-1}\in L^2(\mu_{-1})$, see [C], [S-W] or [S].
\endremark

\heading 4. On the stability of the absence of continuous spectra
 \endheading

Numerous examples given in [A-D], [D], [S] and [R-J-L-S] show that
for a general self-adjoint operator
 the property of not having a continuous part is not stable under
rank one perturbations. For instance [D] contains an example of a
self-adjoint operator $A$ and its cyclic vector $\varphi$ such
that $A$ is pure point but $A+\lambda(\cdot,\varphi)\varphi$ are
singular continuous for all real $\lambda\neq 0$. 

Here we prove that
in order to have $only$ pure point rank one perturbations the operator
must have a very ``thin'' spectrum. We will denote the spectrum of $A$ by $\sigma(A)$.

\proclaim{ Theorem 4.1}
Let $A$ be a self-adjoint  operator. Then the following two 
conditions are equivalent:

1) All self-adjoint rank one perturbations of $A$ are pure point,

2) $\sigma (A)$ is countable.
\endproclaim
\remark{\bf Remark} We assume that $A$ itself is 
included in the set of all its
rank one perturbations. 

In 1) ``all'' can not be replaced with ``almost all'' in any
reasonable sense: it is not difficult
to show that if $\mu$ is a standard singular Cantor measure on $[0;1]$
then $A_\mu+\lambda(\cdot,\varphi)\varphi$ is pure point for a.e. $\lambda\in\Bbb R$ for any cyclic vector $\varphi$.

\endremark

To prove  Theorem 4.1 we will need the following two Lemmas.

\proclaim{Lemma 4.2}
Let $F\subset\Bbb R, \ |F|=0$ be a closed set, $F=\Bbb R\setminus\bigcup_{n=1}^{\infty} I_n$
where $I_n=(x_n;y_n)$ are disjoint open intervals.
Then one can choose two disjoint sets of positive integers $L$ and $M$
such that for any 
$y\in F\setminus \{x_1,x_2,x_3,...\}$
$$\int\limits_ {(y;y+1)\cap \bigcup\limits_{n\in L}I_n}\frac{dx}{x-y}
=\int\limits_ {(y;y+1)\cap \bigcup\limits_{n\in M}I_n}\frac{dx}{x-y}=\infty.\tag4.1$$ and for any
$y\in F\setminus \{y_1,y_2,y_3,...\}$
$$\int\limits_ {(y-1;y)\cap \bigcup\limits_{n\in L}I_n}\frac{dx}{y-x}
=\int\limits_ {(y-1;y)\cap \bigcup\limits_{n\in M}I_n}\frac{dx}{y-x}=\infty.\tag4.2$$
\endproclaim
\demo{Proof}
 Since $|F|=0$ we can choose 
$\{I_{n^{1}_l}\}_{l=1}^\infty$ in such a way
that 
$$\int\limits_ {(y;y+1)\cap 
\bigcup\limits_{l=1}^\infty I_{n^1_l}}\frac 1{x-y}dx>1$$
for any $y\in F\setminus \{x_1,x_2,x_3,...\},$
$$\int\limits_ {(y-1;y)\cap 
\bigcup\limits_{l=1}^\infty I_{n^1_l}}\frac 1{y-x}dx>1$$
for any 
$y\in F\setminus \{y_1,y_2,y_3,...\}$
 and the set of the endpoints
of the intervals $I_{n_l^1}$ has no cluster points except $\pm\infty$. 
Indeed, for each $n\in\Bbb N$ we can choose a finite number of intervals
covering more than a half (in measure) of the interval $[n;n+1]$ . Then we can put
$\{I_{n^1_l}\}$ to be the set of all chosen intervals for all $n\in\Bbb N$.
 
Since the set of the endpoints of $I_{n^1_l}$ has no cluster points in $\Bbb R$, the set
$(0;1)\setminus
\bigcup_{l=1}^\infty I_{n^1_l}$  
is a union of disjoint closed intervals. Therefore in a similar way
we can choose  intervals $I_{n^2_l}
\subset
\Bbb R\setminus
\bigcup_{l=1}^\infty I_{n^1_l}$ 
 such that 
$$\int\limits_ {(y;y+1)\cap 
\bigcup\limits_{l=1}^\infty I_{n^2_l}}\frac 1{x-y}dx>1$$
for any $y\in F\setminus \{x_1,x_2,x_3,...\},$
$$\int\limits_ {(y-1;y)\cap 
\bigcup\limits_{l=1}^\infty I_{n^2_l}}\frac 1{y-x}dx>1$$
for any 
$y\in F\setminus \{y_1,y_2,y_3,...\}$
 and the set of the endpoints of the intervals $I_l^2$ has no cluster points except possibly the ends of the intervals
$I_l^1$ and $\pm\infty$. Proceeding in this way we will obtain
$
\{I_{n^{k}_l}\}_{l=1}^\infty$ for $k=1,2,3...$.
After that we can put 
$L=\{n^{2l-1}_k \}_{l,k\in\Bbb N}$ and $M=\{n^{2l}_k\}_{l,k\in\Bbb N}$.
$\blacktriangle$\enddemo

\proclaim{Lemma 4.3}Let $F\subset (0;1)$ be a closed set, $|F|=0$, $F=(0;1)\setminus\bigcup I_n$
where $I_n=(x_n;y_n)$ are disjoint open intervals. 
Let $\{z_n\}_{n>0}$ be a sequence of points such that $z_n\in I_n$ and $y_n-z_n<|I_n|^2$.
Let $M$ be as in the previous Lemma.
Then there exists $N\subset M$ such that for all $y\in F$

i) if $y=y_n$ for some $n\in\Bbb N$ then
$$\int\limits_{(y;y+1)\cap\bigcup\limits_{n\in N}I_n}\frac {dx}{x-y}=+\infty.\tag4.3$$

ii) if $y\neq x_n,y_n$ for any $n$
and $E=\bigcup\limits_{n\in N}I_n\cup\bigcup\limits_{n\in\Bbb N}[z_n;y_n]$, then
$$p.\ v.\ \int\limits_{(-1;1)}\chi_E(y+x)\frac{dx}{x}<+\infty.\tag4.4$$

\endproclaim
\demo{\bf Proof}
 WLOG we can  assume that the intervals $I_n$ are 
enumerated in such a way that $|I_1|\geq|I_2|\geq...$. 

Let $\delta_l $ be a small positive constant (the exact choice of $\delta_l$ will be made later).
For each $l\in\Bbb N$  choose a sequence of intervals $\{I_{n_k^l}\}_{k>0}$ 
such that for any $k\in\Bbb N$  $\ \ n_k^l\in M$, $I_{n^l_k}\in (y_l;y_l+\delta_l)$,
$\ \ y_{n^1_1}>x_{n^1_1}>y_{n^1_{2}}>x_{n^1_2}>...$ and
$$\int\limits_{(y_l;y_l+1)\cap\bigcup\limits _{k>0}I_k^l}\frac{dx}{x-y_l}=\infty.$$
Put $N=\{n^l_k |\  l,k\in\Bbb N\}$. 
Then condition $i)$ is satisfied.  
We will show that if constants $\delta_l$
are chosen small enough then $N$ also satisfies $ii)$.

Denote $J_l=\bigcup\limits_{k>0}I_{n^l_k}$. Choose constants $\delta_l$ 
to satisfy both 
$$\delta_l<|I_l|^2\tag4.5$$
and
$$\delta_l<\frac1{4}(y_l-z_l)\dist(y_l;(y_l;\infty)\cap\bigcup\limits_{k<l} J_k)\tag4.6$$
for every $l\in\Bbb N$ 
(notice that the distance in (4.6) is always nonzero by the construction).
Let now $y\in F, y\neq x_n,y_n$ for any $n>0$. Consider the sets of integers $M'=\{l|x_l>y\},\ M''=
\{l|\ y_l+\delta_l<y\}$ and $M'''=\{l|\ y_l<y,\ y_l+\delta_l>y\}$. Then $\Bbb N=M'\cup M''\cup M'''$. Put 
$$J'=\bigcup\limits_{l\in M'}J_l,\ \ J''=\bigcup\limits_{l\in M''}J_l \ \ {\text { and}}\ \  J'''=\bigcup\limits_
{l\in M'''}J_l.$$ 
Then, since $J_l\subset (y_l;y_l+\delta_l)$,
$$\left|\int\limits_{J'}\frac{dx}{x-y}\right|<\infty$$
by (4.5). Also obviously
$$\int\limits_{J''}\frac{dx}{x-y}<0<\infty$$
since $J''\subset(y-1;y)$.

Let $\{n_k\}_{k\in\Bbb N}$ be some enumeraton of $M'''$ such that $ n_1<n_2<n_3<...$ . 
Put $\epsilon_k=y-y_{n_k}+\delta_{n_k}$ and consider
$$C_k=\int\limits_{(y-1;y+1)\setminus (y-\epsilon_k;y+\epsilon_k)}\chi_E(x)\frac{dx}{x-y}$$
(notice that $0<\epsilon_k<2\delta_{n_k}\to 0$ as $k\to\infty$).
Then by (4.6) $\delta_{n_k}<\frac14(y_{n_k}-z_{n_k})$ and $\delta_{n_k}<\frac14\dist (y_{n_k};J_l)$ for all $k\in\Bbb N$.
Hence for any
$k\in\Bbb N$ we have that 
$$\dist(y;(y;\infty)\cap\bigcup\limits_{l<k}J_{n_l})>\dist(y_{n_k};(y_{n_k};\infty)\cap\bigcup\limits_{l<k}J_{n_l})-\delta_{n_k}>
$$$$>\frac34\dist(y_{n_k};(y_{n_k};\infty)\cap\bigcup\limits_{l<k}J_{n_l})
>3\frac{\delta_{n_k}}{(y_{n_k}-z_{n_k})}.
$$
Since $|\bigcup\limits_{l<k} J_{n_l}|<1$ and $\delta_{n_k}<\frac14(y_{n_k}-z_{n_k})$, 
$$\frac{|(y;\infty)\cap\bigcup\limits_{l<k} J_{n_l}|}{\dist(y;(y;\infty)\cap\bigcup\limits_{l<k}J_{n_l})}
<\frac{(y_{n_k}-\delta_{n_k})-z_{n_k}}{2\delta_{n_k}}
<\frac{|(z_{n_k};y_{n_k}-\delta_{n_k})|}{\dist\left(y;(z_{n_k};y_{n_k}-\delta_{n_k})\right)}.$$
Therefore
$$\left|\int\limits_{z_{n_k}}^{y_{n_k}-\delta_{n_k}}\frac{dx}{x-y}\right|>
\left|\int\limits_{(y;\infty)\cap\bigcup\limits_{l<k}J_{n_l}}\frac{dx}{x-y}\right|.\tag4.7  $$
Since $y_n-z_n<|I_n|^2$, 
$$\int\limits_
{(y;\infty)\cap\bigcup\limits_{n\in\Bbb N}
(z_n;y_n)}\frac{dx}{x-y}<C<+\infty$$
and by (4.7) 
$$\left|\int\limits_{z_{n_k}}^{y_{n_k}-\delta_{n_k}}\frac{dx}{x-y}\right|>
\left|\int\limits_{(y;\infty)\cap
\left[\bigcup\limits_{n\in\Bbb N}(z_n;y_n)\cup\bigcup\limits_{l<k}J_{n_l}\right]}
\frac{dx}{x-y}\right|-C.\tag4.8  $$
Since $\delta_{n_l}<\frac14(y_{n_l}-z_{n_l})$ for all $l$ and $|I_{n_1}|\geq|I_{n_2}|\geq...$, we have that 
$J_{n_l}\subset(y-\epsilon_k;y+\epsilon_k)$ for all $l\geq k$. 
Thus
$$\left\{\frac{\chi_E(x)}{x-y}>0\right\}\setminus (y-\epsilon_k;y+\epsilon_k)=(y;\infty)\cap
\left[\bigcup\limits_{n\in\Bbb N}(z_n;y_n)\cup\bigcup\limits_{l<k}J_{n_l}\right].$$ Since $$\frac{\chi_E(x)}{x-y}<0$$ on $(z_{n_k};y_{n_k}-\delta_{n_k})$,
(4.8) implies  $C_k<C$ for any $k\in\Bbb N$. Hence 
$$\liminf_{\epsilon\to 0}\int\limits_{(-1;1)
\setminus(-\epsilon;\epsilon)}\chi_E(x)\frac{dx}x<+\infty$$
 and condition (ii) is
satisfied. $\blacktriangle$\enddemo

We will need the following 

\definition{Definition} We will say that two disjoint sets of real
 numbers
$\Cal A$ and $\Cal B$
are well-mixed if they satisfy the following conditions:

1) for any two points $x,y\in\Cal A$ each of the sets $(x;y)$ and $\Bbb R\setminus [x;y]$ 
contains at least one point
from $\Cal B$. 

2) for any two points $x,y\in\Cal B $ each of the sets $(x;y)$ and $\Bbb R\setminus [x;y]$ 
contains at least one point
from $\Cal A$. 
\enddefinition

\demo{Proof of Theorem 4.1}
$2)\Rightarrow 1)$. As follows from the Weyl-von Neumann theorem on the
stability of spectra, the essential spectrum is stable under rank one perturbations.

$1)\Rightarrow 2)$. 
Suppose $\sigma(A)$ is uncountable. 
We will show that there exists a self-adjoint rank one perturbation
with nontrivial continuous part.

If $A$ is not cyclic 
then 
we can always choose a cyclic
subspace such that the restriction of $A$ on this subspace has an uncountable spectrum. 
Indeed, $A$ is unitarily equivalent to the multiplication by $z$ in
a direct integral of Hilbert spaces $H_\xi$ 
$$H=\int\oplus H_\xi d\mu(\xi)$$
where $\mu$ is  a scalar measure with uncountable support.
If we choose a one-dimensional subspace $l_\xi$ in each of these Hilbert spaces and consider the direct
integral 
$$L=\int\oplus l_\xi d\mu(\xi)$$ 
then the restriction of multiplication by $z$
on $L$ will be cyclic with the spectral measure $\mu$. Since $\supp\mu$ is uncountable, 
the spectrum of such a restriction will be uncountable. 

Notice that if some  
rank one perturbation of this restriction has a nontrivial continuous part then
the corresponding  rank one perturbation of the whole operator also has a nontrivial continuous part.
Therefore, WLOG we can assume that $A$ is cyclic.

 Denote by $\Cal A$
the set of all eigenvalues of $A$. WLOG $\Cal A\subset [0;1]$. Then 
there exists a closed uncountable set $F\subset(0;1), m(F)=0$
such that for any 
$x\in F$ and any $\epsilon>0$ both sets  $ (x -\epsilon;x)\cap \a$
and   $ (x;x +\epsilon)\cap \a$ are nonempty.

Let $I_1=(x_1;y_1),I_2=(x_2;y_2),...$ 
be disjoint open intervals such that $F=(0;1)\setminus\bigcup                        
I_n$. We always can assume that $F$ contains no isolated points i.e.
$x_i\neq y_j$ for any $i,j$.

First, using Lemmas 4.2 and 4.3 we will construct a 
 phase shift $u_0$ of a pair $(\mu_0;\nu_0)$ such that $\mu_0=\sum c_n\delta_{a_n}$
for some  $\ a_n\in\a, c_n>0$ and $\nu_0$ is continuous with $\supp\nu_0\subset F$.

By Lemma 4.2 we can choose two sets of integers $L$ and $M$ satisfying (4.1) and (4.2).
For each $n\in\Bbb N$ choose $z_n\in I_n\cap\Cal A$ such that $y_n-z_n<|I_n|^2$. 
By Lemma 4.3
we can choose $N\subset M$ satisfying (4.3) and (4.4). Denote $O=\Bbb N\setminus N$.
For each $ k\in O$ put $a_k=z_k$. 
For each $ k\in N$ choose $a_k\in I_k\cap\a$ such that 
$a_k-x_k<|I_k|^2$. After that  put
$u_0=\pi$ on $\bigcup(a_k;y_k)$ and $u_0=0$ elsewhere.
Then $u_0$ has compact support.

Consider the pair of measures $(\mu_0;\nu_0)$ such that $u_0$ is its phase shift.
By Lemma 3.2 $\mu_0$ and $\nu_0$ are singular.
To prove that $\mu_0=\sum a_n\delta_{a_n}$ we need to show that $\mu_0(F)=0$.

To show that $\mu_0(F)=0$ let us notice that 
if $y\in F$ then
$$0\geq\int\limits_{(y-1;y)\cap\bigcup\limits_{k\in N}(x_k;a_k)}\frac{dx}{x-y}\geq -C_1>-\infty$$
because
$a_k-x_k<|I_k|^2$ for any $k\in N$. Thus
$$\int\limits_{\{u_0=\pi\}}\frac{dx}{x-y}<\int\limits_{\bigcup\limits_{k\in\Bbb N} (z_k;y_k)\cup\bigcup\limits_{k\in N}I_k}\frac{dx}{x-y}+C_1.
\tag4.9$$
If $y\neq y_i,x_i$ for any $i\in\Bbb N$ then
(4.9) and  Lemma 4.3 imply that there exist positive $\epsilon_1,\epsilon_2,...$ such that
$\epsilon_n\to 0$ as $n\to\infty$ and 
$$\int\limits_{\{u_0=\pi\}\setminus(y-\epsilon_n;y+\epsilon_n)}\frac{dx}{x-y}<C_2<\infty.$$
for any $n$. Hence,
$$p.\ \ v.\ \ \int\limits_{\Bbb R}\frac{u_0(y+x)dx}{x}<\infty.$$
Therefore  
by Lemma 3.3 $\ \ \mu_0(F\setminus\{x_1,y_1,x_2,y_2,...\})=0$.
Also, $\mu_0(x_k)=\mu_0(y_k)=0$ for all $k$ because $u_0(x)\to0$ as $x\to y_k+$
and $u_0(x)\to\pi$ as $x\to x_k-$ which means that condition (3.8) is not satisfied.
Hence, $\mu_0=\sum c_n\delta_{a_n}$ for some positive $c_n$.

For $\nu_0$ we obviously have $\supp\nu_0\subset F$. Let $y\in F, y=y_i$ for some $i>0$.
Then
$$\int\limits_{(y;y+1)\cap\bigcup\limits_{k\in N}I_k}\frac{dx}{x-y}=\infty$$
because $N$ satisfies (4.3).
Since for $k\in N$ we have $a_k-x_k<(y_k-x_k)^2$, this implies that
$$\int\limits_{(y;y+1)\cap\bigcup\limits_{k\in N}(a_k;y_k)}\frac{dx}{x-y}=\infty.$$
Since on $(a_k;y_k)$ $u_0=\pi$, condition (3.9) is not satisfied. Thus $\nu_0(y)=0$.
If $y\in F$ and $y\neq y_i$ for any $i$, then
$$\int\limits_{(y-1;y)\cap\bigcup\limits_{k\in O}I_k}\frac{dx}{y-x}=\infty$$
by (4.2) because $M\subset O$. Since for $k\in O$ we have $y_k-a_k<(y_k-x_k)^2$, this implies
that
$$\int\limits_{(y-1;y)\cap\bigcup\limits_{k\in O}(x_k;a_k)}\frac{dx}{y-x}=\infty.$$
Again, since $u_0=0$ on $(x_k;a_k)$, condition (3.9) does not hold and $\nu_0(y)=0$.
Hence, $\nu_0$ is continuous.

Our next goal is to transform $u_0$ into the phase shift $u$ of a pair $(\mu;\nu)$ such
that $\mu$ is a spectral measure of $A$ and $\nu$ has a nontrivial continuous part
(note that $\mu_0$ does not have point masses at some points of
 $\a$, therefore it is
not  a spectral measure of $A$).

Let $\{b_n\}_{n=1}^\infty$ be some enumeration of the set $\a\setminus
\{a_n\}_{n=1}^\infty$. For each $b_n$ let us choose $c_n$
such that if $b_n\in I_k$ then $c_n\in I_k$,

 $$\ \ |b_n-c_n|<\frac{\dist(b_n;(\Bbb R\setminus I_k)
\cup\{b_1,b_2,...,b_{n-1}\})}{2^{n+1}}\tag4.10$$ and
 for each $n\in\Bbb N$ the sets 
$\Cb_n=\{a_i\}_{i=1}^{\infty}\cup
\{b_1,b_2,...,b_n\}$ and $\Cc_n=\{c_1,c_2,...,c_n\}\cup F$ are well-mixed.

For each $k\in\Bbb N$ define the function $u_k$ in the following way:

1) $|u_k|=0$ or $\pi$ everywhere on $\Bbb R$;

2) $u_k$ is continuous everywhere except $ \Cb_k\cup\Cc_k$;

3) $u_k$ jumps from $0$ to $\pi$ at each point of $\Cb_k$
and from $\pi$ to $0$ at each point of $\{c_1,c_2,...,c_k\}$.

Let $u_k$ be the phase shift of a pair $(\mu_k;\nu_k)$.
Then 
$$\mu_k=\sum _{n=1}^\infty\alpha^k_n\delta_{a_n} 
+\sum_{n=1}^k \beta^k_n\delta_{b_n}$$ 
for some positive constants
$\alpha^k_i$ and $\beta^k_i$ and 
$$\nu_k=\sum_{n=1}^k\gamma^k_n \delta_{c_n}+f_k\nu_0$$
for some positive constants $\gamma^k_1$ and some positive function 
$f_k\in L^1(\nu_0)$. Condition (4.10) implies that the sequence 
$\{u_k\}$ converges in measure to some function $u$. Since $\supp u_k\subset [-1;2]$
for any $k\in\Bbb N$, $\ \ u $ has compact support.
Let $u$ be the phase shift of a pair $(\mu;\nu)$.
Then $\mu_k\to\mu$ and $\nu_k\to\nu$ in the $*$-weak topology.
 
Since by (4.10)
 $$\int\limits_{\Bbb R}\left|\frac {u_k(x)-u_{k+1}(x)}{x-y}dx\right|<\frac 1{2^{k+1}}\tag 4.11$$
for any $y\in\Cb_k$, 
by part ii) of Lemma 3.7
$$\mu_{k+1}=g_k\mu_k+\beta^{k+1}_{k+1}\delta_{b_{k+1}}$$
 where
$g_k\in L^\infty(
\mu_k)$, 
$$1-\frac1{2^k}<g_k<1+\frac1{2^k}\tag4.12$$ $\mu_k$-a. e. 
Also by (4.11) and part ii) of Lemma 3.7
$$1-\frac1{2^k}<\frac{f_k}{f_{k+1}}<1+\frac1{2^k}\tag4.13$$
$\nu_0$-a. e.

Since 
$$||\mu_k||=
\im \exp \left[\Cal Ku_k(i)\right]
\to \im \exp \left[\Cal Ku(i)\right],$$
$||\mu_k||\to||\mu||$ and  $\ \beta^{k+1}_{k+1}\to 0$ as $k\to\infty$. 
Since  each $\mu_k$ is pure point, that implies, together with (4.12),
that $\mu$ is pure point and that
$$\mu=\sum _{n=1}^\infty\alpha_n\delta_{a_n} 
+\sum_{n=1}^\infty \beta_n\delta_{b_n}$$ where
$$0<\alpha^1_i\prod_{k\in\Bbb N} g_k(a_i)=\alpha_i<\infty$$ and 
$$0<\beta^i_i\prod_{k>i} g_k(b_i)=\beta_i<\infty.$$ Thus $A$ is unitarily
equivalent to $ A_\mu$.
Also, since $\nu=\eta+\sigma$ where $\eta$ is some positive measure
and $\sigma$ is a $*$-weak limit of the sequence $\{f_k\nu_0\}$,
(4.13) implies that $\sigma=f\nu_0$ for some $f\in L^1(\nu_0)$, 
$$0<f_1(x)\prod_{k>1}f_k(x)=f(x)<3f_1(x)$$
$\nu_0$-a. e. Thus $A_\mu+(\cdot,1)1$ has a nontrivial continuous part.
$\blacktriangle$\enddemo

\proclaim{Corollary 4.4}
Let $A$ be a self-adjoint operator. Then all trace class perturbations
of $A$ are pure point iff $\sigma(A)$ is countable.
\endproclaim
\demo{Proof} 

The ``if'' part follows from the stability of the essential spectrum under trace class
perturbations.

The ``only if'' part follows from Theorem 4.1.
 
The statement also follows from the Theorem of Carey and Pincus [C-P] on the equivalence
modulo the trace class. $\blacktriangle$ \enddemo 

\heading 5. The problem of two spectra \endheading

\define\rank{\operatorname {rank}}
\definition{Definition}
We will say that two operators $A$ and $B$ are 
equivalent up to a rank one perturbation 
 if 
there exist operators $A'$ and $B'$ acting in the same space
such that $A$ is unitarily
equivalent to $A'$, $B$ is unitarily equivalent to $B'$ and 
$\rank (A'-B')=1$.
\enddefinition

In this section we will give  a partial answer 
to the following question.

\remark{The problem of two spectra}
Let $\mu$ and $\nu$ be two finite Borel measures on $ {\Bbb R}$.
When do there exist two cyclic self-adjoint operators $A$ and $B$
 equivalent up to a rank one perturbation 
such that $\mu$ and $\nu$ are spectral measures
of $A$ and $B$ respectively?
\endremark

\definition{Definition} If such $A$ and $B$ exist, we will
say that $\mu$ and $\nu$ solve  the problem of two 
spectra (PTS).
\enddefinition

 We will say that $\mu$ is equivalent to $\nu$ 
and write $\mu\sim\nu$ if there exists
$f\in L^1(\mu)$ such that $f>0\ \ \  \mu$-a. e. and $f\mu=\nu$. 
We will 
say that  pairs of measures  $(\mu,\nu)$ and  $(\mu',\nu')$ are equivalent,
 $(\mu,\nu)\sim(\mu',\nu')$, if  $\mu\sim\mu'$ and  $\nu\sim\nu'$.

In terms of the families $M_\phi$, we can say that $\mu$
and $\nu$  solve PTS if and only if there exists $\phi$
such that $(\mu,\nu)\sim(\mu',\nu')$ for some $\mu',\nu'\in\Cal M_\phi$.
In terms of the phase shift, $\mu$
and $\nu$  solve PTS if and only if there is an equivalent pair possessing
a phase shift i. e. satisfying (3.3) for some real  
function $u\in L^\infty(\Bbb R), \ ||u||_\infty\leq\pi$ and $\lambda\in\Bbb R$.

In this section we will discuss the case of pure point 
measures $\mu$ and $\nu$.

The first result in this direction is the following theorem, proved by
Gelfand and Levitan.

\proclaim{Theorem 5.1}
Let $\Cal A=\{a_n\}_{n=1}^\infty$ and $\Cal B=\{b_n\}_{n=1}^\infty$ be two disjoint  
sequences
of real numbers,  $\lim_{n\to\infty}a_n=\lim_{n\to\infty}b_n=c\in\Bbb R$ .
There exist unitarily equivalent up to a rank one perturbation
cyclic self-adjoint operators $A$ and $B$
 such that $\sigma(A)=\Clos\Cal A$ and  
$\sigma(B)=\Clos\Cal B$
if and only if the sets $\{a_1,a_2,...\}$ and $\{b_1,b_2,...\}$ are well-mixed.

If such $A$ and $B$ exist then they are
 unique up to a unitary equivalence.
\endproclaim

\demo{Proof} 
Define function $u$ on $\Bbb R$ to satisfy the following conditions:

1) $u$ is continuous and $|u|$ is equal to
$\pi$ or 0  everywhere on 
$\Bbb R\setminus (\{a_n\}\cup\{b_n\})$,

2) at the points $\{a_n\}$ 
$u$  jumps from $0$ to $\pi$,

3) at the points
$\{b_n\}$ $u$  jumps from $\pi$ to $0$.

WLOG 
$\inf_{n\in\Bbb N}a_n\leq\inf_{n\in\Bbb N}$ and  $\sup_{n\in\Bbb N}a_n\leq \sup_{n\in\Bbb N}b_n$.
Then $u$ has compact support. 
As we discussed before, there are unique measures $\mu,\nu\in\Cal M_+(\Bbb R)$ such that
$u$, $\mu$ and $\nu$ satisfy formula (3.3) with $\lambda=1$.
Put $A=A_\mu$ and $B=A_\nu$. Then $B$ is unitarily equivalent to $ A+(\cdot,1)1$. 
By Lemma 3.3 i) 
$\supp\mu\subset\Clos\Cal A$ and $\supp\nu\subset\Clos\Cal B$. 
Also by Lemma 3.8
the measures $\mu$ and $\nu$ have point masses at points $a_i$ and $b_i$ respectively, therefore
$\supp\mu=\Clos\Cal A$ and $\supp\nu=\Clos\Cal B$.
Thus the operators
$A$ and $B$ satisfy the condition of the Theorem.
Since the sequences are well-mixed and have only one cluster point
there exists a unique function $u$ satisfying properties 1)-3).
Therefore $A$ and $B$ are unique up to a unitary equivalence.
$\blacktriangle$\enddemo

\remark{\bf Remark} 
If
$\mu=\sum \alpha_n\delta_{a_n}$
and $\nu=\sum \beta_n\delta_{b_n}$, where $\alpha_i,\beta_i>0$,  solve PTS then the sets
$\Cal A=\{a_1,a_2,...\}$ and $\Cal B=\{b_1,b_2,...\}$ are well-mixed. To prove it one
can, for instance, notice that
if $\Cal A$ and $\Cal B$ are not well-mixed then the function
$u$ satisfying properties 1)-3) from the above proof does not exist.

Note that Theorem 5.1  provides us with operators
such that $\mu$ and $\nu$  are only {\it absolutely  continuous 
} with  respect to spectral measures of $A$ and $A+(\cdot,\phi)\phi$ but not
necessarily {\it equivalent} to them. Therefore it does not
imply that $\mu$ and $\nu$   solve PTS. 
As shown in Example 5.2 below, 
 the condition that $\Cal A$ and $\Cal B$ are well-mixed
is not sufficient for $\mu$ and $\nu$
to   solve PTS.
\endremark

\proclaim{Example 5.2}
Put $ a_n=(-1)^n/2^n$ for $n=1,2,3,...$, $b_1=-1,\ b_n=a_{n-1}+(-1)^n/4^n$ 
for $n=2,3,...$. 
Then the sets $\{a_1,a_2,...\}$ and $\{b_1,b_2,...\}$ are well-mixed.

Let $\mu$ and $\nu$ be defined as in the last remark.
Suppose that $\mu$ and $\nu$  solve PTS. Then there exist a
 pair  $(\mu,\nu)\sim(\mu',\nu')$ possessing a phase shift $u$.
Such $u$ must be continuous on $\Bbb R\setminus
(\{a_n\}\cup\{b_n\})$, jump from $0$ to $\pi$ at any $a_n$
and jump from $\pi$ to $0$ at any $b_n$.
Thus
$ u=0$ on $(1/2^{2n}-1/4^{2n+1};1/2^{2n})$ 
and on $(-1/2^{2n-1}+1/4^{2n};-1/2^{2n+1})$ for $n=1,2,...$
and $u=\pi$ 
on the rest of $\Bbb R$. 

But then  Lemma 3.8 implies that $\mu(0)>0$. Since $0\neq a_n$ for any $n\in\Bbb N$ we have a contradiction.
\endproclaim

Hence in Theorem 5.1 one can not  prescribe arbitrarily whether or not spectral
measures of $A$ and $B$ have a point mass at $c$.

If the sequences $\{a_1,a_2,...\}$ and $\{b_1,b_2,...\}$  are finite, disjoint 
and well-mixed then $\mu$ and $\nu$ obviously  solve PTS. More interesting example
is provided by the following theorem.

\proclaim{Theorem 5.3 ([A3])}
Let $\Cal A$ and $\Cal B$ be two disjoint countable sets
on the unit circle $\Bbb T$. Suppose that 
${\Clos\Cal A}=\Clos{\Cal B}=\Bbb T$.
Let $\{a_n\}_{n=1}^\infty$ and $\{b_n\}_{n=1}^\infty$ be some enumerations of
$\Cal A$ and $\Cal B$ respectively.
Then there exist  sequences of positive real numbers $\{\alpha_n\}$
and $\{\beta_n\}$ such that the measures $\sum \alpha_n\delta_{a_n}$
and $\sum \beta_n\delta_{b_n}$ belong to the same family $M_\varphi$
\endproclaim

In terms of operator theory Theorem 5.3 means that any two cyclic pure point self-adjoint
operators whose spectrum is equal to $\Bbb R$ are unitarily 
equivalent up a rank one perturbation.
Note that the  condition ${\Clos\Cal A}=\Clos{\Cal B}=\Bbb T$ automatically implies that
$\Cal A$ and $\Cal B$ are well-mixed.
However the following 
example shows 
 that in general the
condition $\Clos{\Cal A}=\Clos{\Cal B}$ does not imply that 
that the corresponding measures  solve PTS even if  $\Cal A$ and $\Cal B$
are well-mixed.

\proclaim{Example 5.4}
Let $C$ be the standard Cantor null set on the unit interval
$[0;1]$, $C=[0;1]\setminus \bigcup I_n$ where $I_n=(x_n;y_n)$ are disjoint
open intervals. 
 Let $\Cal A$ and $\Cal B$ be two disjoint 
well-mixed countable sets of points of $C$ such that 
$$ 0, x_1, x_2,..\in \Cal A, 1\in\Cal B.\tag 5.1$$ Then 
$\Clos{\Cal A}=\Clos{\Cal B}=C$.
Let $\{a_n\}_{n=1}^\infty$ and $\{b_n\}_{n=1}^\infty$ be some enumerations of
$\Cal A$ and $\Cal B$ respectively.
Define $\mu=\sum \frac 1{2^n}\delta_{a_n}$
and $\nu=\sum \frac 1{2^n}\delta_{b_n}$. Suppose that $\mu$ and $\nu$
 solve PTS. Then there exists the phase shift $u$ of the pair
$(\mu';\nu')$ where $(\mu,\nu)\sim(\mu',\nu')$. 
Function $u$ must be constant on each $I_n$.
Condition (5.1) implies that $u=0$ a. e. on $\Bbb R$. This means that $\mu'$ and $\nu'$
are trivial   and we obtain a contradiction.
\endproclaim

It turns out that this example is, in a certain sense, typical:
the following Theorem shows that  such $\mu$
and $\nu$  solve PTS for any two well-mixed dense in $K$ sequences
 if and only if the set $K$ is not too ``porous''. 

\proclaim{Theorem 5.5}
Let $K\subset\Bbb R$ be a closed set.
Denote by $I_1,I_2,...$ the
disjoint open intervals $I_k=(x_k;y_k)$ such that
$K=\Bbb R\setminus\bigcup I_n$. 

Then the following two conditions are equivalent:

i)Any two cyclic self-adjoint pure point operators
$A$ and $B$ such that 
 the sets of their eigenvalues are well-mixed
and $\sigma(A)=\sigma(B)=K$ are unitarily equivalent up to a  rank one perturbation.

ii) If $y\in 
K\setminus\{x_1,y_1,x_2,y_2,...\}$  then 
$$\int\limits_{\bigcup\limits_{n\in \Bbb N} I_n} \frac{dx}{|y-x|}<\infty;\tag 5.2$$
if $y=x_k$ or $y=y_k$ for some $k\in \Bbb N$ then
$$\int\limits_{\bigcup\limits_{n\in \Bbb N, n\neq k} I_n} \frac{dx}{|y-x|}<\infty.\tag 5.3$$
\endproclaim
\example{\bf Example 5.6}
Condition ii) fails for $K=C$ for any standard Cantor set $C$
(even if $|C|>0$). 
Condition ii) obviously holds if $\partial K$ is a finite set.
We will show that there exist more complicated examples of such $K$, namely
that for any nowhere dense set $L\subset\Bbb R$ there exists a set $K$
such that $K$ satisfies ii) and $\partial K\supset L$.

Denote by $J_i=(x_i;y_i)$ the disjoint open intervals such that $\Clos L=\Bbb R\setminus\bigcup J_n$ .
On each $J_i$ choose a sequence of intervals $\{\Delta_i^k\}_{-\infty}^\infty, \Delta_i^k=(x_i^k;y_i^k)$
such that $x_i^k\to x_i$ as $k\to-\infty$, $x_i^k\to y_i$ as $k\to\infty$ and
$$\int\limits_{\cup_{k}\Delta_i^k}\frac{dx}{|x-y|}<\frac1{2^i}$$
for all $y\in\Bbb R\setminus J_i$. Then we can put $K=\Bbb R\setminus \bigcup_{i,k}\Delta_i^k$.
 
Similarly, one can show that for any closed $L\subset \Bbb R$ there exists a closed $K\supset L$
satisfying ii) and such that $\partial K\supset\partial L$. 

\endexample

\demo{Proof of Theorem 5.5}
ii)$\Rightarrow$i)
WLOG $K\subset(0;1)$.
Let $\Cal A$ and $\Cal B$
be the sets of eigenvalues of $A$ and $B$ respectively.
Let $\{a_n\}_{n=1}^\infty$ and $\{b_n\}_{n=1}^\infty$ be some enumerations of
$\Cal A$ and $\Cal B$ respectively.

a) Let us first consider the case when the endpoints of the intervals $I_k$
 belong to neither $\Cal A$ nor $\Cal B$.

Let us find reenumerations $\{a_{n_k}\}$
and $\{b_{m_k}\}$ of
 the sequences $\Cal A=\{a_n\}$ and $\Cal B=\{b_m\}$ 
such that
for any $k\in\Bbb N$

1) 
the sets
 $\{a_{n_1},...,a_{n_{k}}\}$ and $\{b_{m_1},
...,b_{m_{k}}\}$ are well-mixed and 

2) the distance  $|a_{n_k}-b_{m_k}|$ is at least $2^{k+1}$ times less than
the distance between the $a_{n_k}$ and any of the points
 $a_{n_1},...,a_{n_{k-1}},b_{m_1},
...,b_{m_{k-1}}$.

To find such reenumerations we can act
in the following way.

WLOG $a_1>b_1$.
Choose $n_1=1$ and $m_1=1$. After that for $k=2,3,...$  consequently do the following.
Choose a point $x$ in  $\Cal A\cup\Cal B\setminus 
 \{a_{n_1},...,a_{n_{k-1}},b_{m_1},
...,b_{m_{k-1}}\}$ with minimal index in the initial enumeration. Suppose $x\in \Cal A$. Then put $a_{n_k}=x$
and choose $b_{m_k}\in\Cal B$ to satisfy conditions 1) and 2) above. Similarly,
if $x\in\Cal B$ then  put $b_{m_k}=x$
and choose $a_{n_k}\in\Cal A$ to satisfy conditions 1) and 2) above. 

Note that it is always possible to choose $b_{m_{k}}$ or  $a_{n_{k}}$ to satisfy 2) 
 because $a_{n_k},b_{m_{k}}\in K\setminus \{x_1,y_1,x_2,y_2,...\}$.

Now for each $k$ define a function $u_k$ on $\Bbb R$ in such a way
that $|u_k|\in\{0,\pi\}$ on $\Bbb R$, $u_k$ is
continuous everywhere except $a_{n_1},a_{n_2},...a_{n_k},b_{n_1},b_{n_2},...,b_{n_k}$
, $u_k$ jumps from $0$ to $\pi$ at each $a_{n_i}, 0<i\leq k$
 and $u_k$ jumps
from $\pi$ to $0$ at each $b_{n_i}, 0<i\leq k$ (this construction is possible
because the sequences are well-mixed).
Since $a_1>b_1$,
$u_k$ has a compact support. Therefore $u_k$ 
is the phase shift of a pair of measures $(\mu_k;\nu_k)$.
It is easy to see that
$\mu_k=\sum_{i=1}^k \alpha^k_i\delta_{a_{n_i}}$ 
and
$\nu_k=\sum_{i=1}^k \beta^k_i\delta_{b_{m_i}}$ 
for some positive constants 
$\alpha^k_1,...,\alpha^k_k,\beta^k_1,...,\beta^k_k$.

Condition 2)  implies that for any $l>k$ 
$|\{u_k\neq u_l\}|\leq
\sum_{i=k}^l 1/2^i$.
Thus the sequence $\{u_k\}$ converges in measure to some function $u$.
That implies that the sequences of measures $\{\mu_k\}$ and $\{\nu_k\}$  converge in the 
$*$-week topology to some measures $\mu$ and $\nu$ (respectively) such that $u$ is the phase shift
of the pair $(\mu;\nu)$.
To complete the proof we need to show that
$\mu=\sum \alpha_n\delta_{a_n}$
and 
$\nu=\sum \beta_n\delta_{b_n}$
 for some positive constants $\alpha_i$ and $\beta_i$.

Since 
$$||\mu_k||=
\im \exp \left[\Cal Ku_k(i)\right]
\to \im \exp \left[\Cal Ku(i)\right],$$
we have that $$||\mu_k||-||\mu_{l}||\to 0$$
as $k,l\to\infty$. Thus 
$$
\sum_{i=k}^l\alpha^{i}_{l}\leq\left|\  
||\mu_k||-||\mu_l||\ \right|+\sum_{i=1}^k|\alpha^i_k-\alpha^i_l|
\leq\left|||\mu_k||-||\mu_l||\right|+k/2^{k-1}\to 0\tag5.4$$
as $k,l\to\infty$. 
But condition 2) implies that 
$$\int\limits_{\{u_k\neq u_{k+1}\}}\frac{dx}{|x-a_i|}<\frac1{2^k}$$
for any $k$ and any $i\leq k$. Hence, by part ii) of Lemma 3.7 
$$|\alpha^k_i/\alpha^{k+1}_i-1|<1/2^{k-1}\tag5.5$$
for any $k$ and any $ i\leq k$. 
Since $\supp u_k\subset (0;1)$,
$\alpha^{k+1}_i\leq ||\mu_k||\leq 1$. Hence by (5.5)
$$|\alpha_i^k-\alpha_i^{k+1}|<1/2^{k-1}.\tag5.6$$

Together (5.4) and  (5.6) imply that
$$||\mu_k-\mu_l||\to 0$$ 
as $k,l\to\infty$. 
Put 
$\alpha_i=\lim_{k\to\infty}\alpha_i^k$, 
and
$\mu'=\sum_{i=1}^\infty\alpha_i\delta_{a_{n_i}}$.
Since the sequence $\{\mu_k\}$ converges in norm, it must converge to $\mu'$.
But since $\mu$ is the $*$-weak limit of $\{\mu_k\}$, 
$\mu=\mu'=\sum_{i=1}^\infty\alpha_i\delta_{a_{n_i}}$. It is left to notice
that (5.5) implies that  
$$\alpha_i=\alpha_i^i\prod_{l>i}\frac {\alpha_i^l}{\alpha_i^{l+1}}>0$$
for each $i$.
In the same way we can prove that $\nu=\sum \beta_i\delta_{b_i}$ for some positive $\beta_i$.

b) In the general case, 
denote 
 $N=\{n| \partial I_n\cap K\neq\emptyset\}$.
Consider disjoint countable sets 
 $\Cal A'=\Cal A\cup\Cal C$ and $\Cal B'=\Cal B\cup\Cal D$ where
$\Cal C=\{c_n\}\subset[\cup_{n\in N}I_n]\cap (0;1)$ 
 and $\Cal D=\{d_n\}\subset [\cup_{n\in N}I_n]\cap (0;1)$ satisfy:

1) $\Cal C,\Cal D$  are dense in 
$[\cup_{n\in N}I_n]\cap (0;1)$ 
and 

2) if for some $n\in N$ one of the endpoints $x$ of the interval $I_n$ belongs to $\Cal A\ (\Cal B)$
but the other endpoint $y$ of $I_n$ is not in $\Cal A\cap\Cal B$ then $y\in \Cal A'\ (\Cal B')$
(i. e. both $x$ and $y$ are in the same set $\Cal A'$ or $\Cal B'$).

(Notice that since $\Clos\Cal A=\Clos\Cal B=K$ and $\Cal A\cap\Cal B=\emptyset$, $K$ does not have isolated
points. Hence $x_k\neq y_l$ for any $k, l \in \Bbb N$.)

Then the set $\Cal A'\cup\Cal B'$ contains all the endpoints of the intervals $I_n,\ n\in\Bbb N$.

Then $K'=K\cup\bigcup_{n\in N}I_n,\  \Cal A'$ and $\Cal B'$ satisfy the conditions of part a) and we
can construct a phase shift $v$ corresponding to measures $\mu'=\sum\alpha_i\delta_{a_i}+\sum\sigma_i\delta_{c_i}$
and $\nu'=\sum\beta_i\delta_{b_i}+\sum\eta_i\delta_{d_i}$ for some positive
$\alpha_i,\sigma_i,\beta_i,\eta_i$. 
Define $N_\pi=\{n| x_n\in\Cal A\  or\  y_n\in\Cal B\}$ and  
$N_0=\{n| x_n\in\Cal B\  or\  y_n\in\Cal A\}$. Then 
 $N_\pi\cap N_0=\emptyset$ because the sets $\Cal A$ and $\Cal B$ are well-mixed and
$N=N_\pi\cup N_0$.
Let  function $u$ be such that
$u=v$ outside $\cup_{n\in N}I_n$, 
$u=\pi$ on each $I_n, n\in N_\pi$ 
and $u=0$ on each $I_n, n\in N_0$. Let
$\mu$ and $\nu$ be the measures such that $u$ is the phase shift of the pair
$(\mu;\nu)$. Then ii) together with Lemma 3.7 imply that the restriction of $\mu$ on
$K\setminus\{x_1,y_1,x_2,y_2,...\}$ is absolutely continuous with respect to
$\mu'$. Also, obviously $\mu(\cup_{n\in N} I_n)=0$ because $u$ is constant on each
$I_n, \ n\in N$. It is left to check the endpoints of the
intervals $I_n,  n\in N$. 

If for some $n\in N_\pi$ $x_n\in \Cal A$, then $\mu'$ has a point mass
at $x_n$. Thus $v$ satisfies the condition (3.8) at the point $y=x_n$. Since $u=\pi$
on $I_n$, by (3.8) and ii) we have
$$\int\limits_{\{u\neq v\}}\frac{dx}{|x-x_n|}<\infty.$$
Hence by (3.8) $\mu$ also has a point mass at $x_n$.

If for some $n\in N_\pi$ $x_n$ belongs to neither $\Cal A$ nor $\Cal B$
 then $\nu'$ has a point mass at $x_n$. Thus $v$ satisfies
(3.9) at $y=x_n$. Since $u=\pi$ on $I_n$, by (3.8) and ii) we have that
$$\int\limits_{\{u\neq \pi\}}\frac{dx}{|x-x_n|}<\infty.$$
Thus by (3.8) and (3.9)
$\mu(x_n)=0$ and $\nu(x_n)=0$.

Other endpoints of $I_n,\ n\in N$ can be checked in the same way.

Thus
$\mu=\sum \alpha_n\delta_{a_n}$ for some positive $\alpha_n$. Similarly we can show
that $\nu=\sum \beta_n\delta_{b_n}$ for some positive $\beta_n$.

i)$\Rightarrow$ii)

Suppose condition ii) is not satisfied at some $y\in K$. Let $\Cal A$ and $\Cal B$ be disjoint,
dense in $K$ and well-mixed sequences such that $y\in\Cal B$ and $y_i\in\Cal A$ for all $i$ such that
$y_i\neq\infty$.
Then the phase shift $u$ of the 
pair $(\mu;\nu)$ is equal to $0$ on $\cup I_n$.
Thus condition (3.9) is not satisfied at $y$ and 
$\nu(y)=0$. $\blacktriangle$ \enddemo

Here is one more way 
to avoid the situation of Example 5.4.

\proclaim{Theorem 5.6}
Let $K\subset\Bbb R$ be a closed set.
Denote by $I_1=(x_1;y_1),I_2=(x_2;y_2),...$ disjoint open intervals such that $K=\Bbb R\setminus\bigcup I_n$.
Let
$A$ and $B$ be two self-adjoint cyclic pure point operators, $\Cal A$ and $\Cal B$ be the sets
of all eigenvalues of $A$ and $B$ respectively. Suppose that $\sigma (A)=\sigma(B)=K$ and 
$\Cal A\cap\{x_1,y_1,x_2,y_2,...\}=
\Cal B\cap\{x_1,y_1,x_2,y_2,...\}=\emptyset$. Then $A$ and $B$ are equivalent
up to a rank one perturbation.
\endproclaim
\demo{Proof}
Notice that in part a) of the implication ii)$\Rightarrow$i) in the previous
proof we did not use condition ii). $\blacktriangle$
\enddemo

In the rest of this Section we are going to discuss the following 
question.

\remark{\bf Remark}
Let $A$ be a singular self-adjoint cyclic operator, $\varphi$ and $\psi$  its
noncollinear cyclic vectors  ($\varphi\neq c\psi$). Is it possible that
operators $A^\varphi=A+(\cdot,\varphi)\varphi$ and $A^\psi=A+(\cdot,\psi)\psi$
are unitarily equivalent?

If $A$ is  a finite rank operator, then the answer is negative.
Indeed, denote by $\mu$ and $\nu$ the spectral measures of $\varphi$ for
$A$ and $A^\varphi$ respectively; denote by $\mu'$ and $\nu'$ the spectral measures of
$\psi$ for $A$ and $A^\psi$ respectively.
If $\mu$ and $\nu$ are linear combinations of point masses
at points $a_1,a_2,..,a_n$ and $b_1,b_2,..,b_n$ respectively,
$a_1<b_1<a_2<b_2<...<a_n<b_n$, then the phase shift $u$ of the pair
$(\mu;\nu)$ depends only on the sequences
$\{a_n\}$ and $\{b_n\}$: $u=0$ a. e. on $(-\infty;a_1)\cup(b_1;a_2)\cup
...\cup(b_{n-1};a_n)\cup(b_n;\infty)$ and $u=\pi$ elsewhere. 
If $A^\varphi$ is unitarily equivalent to $A^\psi$ then $\nu\sim\nu'$. Also we have that $\mu\sim\mu'$.
Thus measures $\mu'$ and $\nu'$ have point masses at the same points $a_1,a_2,...,a_n, b_1,b_2,...,b_n$.
Hence the phase shift $u'$ of the pair $(\mu',\nu')$ must be equal to $u$ a. e. on 
$\Bbb R$. Since the phase shift  determines the pair of measures uniquely, we have $\mu=\mu'$ and $\nu=\nu'$.
Therefore, $\varphi=c\psi$ and we have a contradiction.
In the same way one can prove that if all eigenvalues of $A$ are isolated, then any two
different (corresponding to noncollinear vectors) rank one perturbations of $A$  can not be
unitarily equivalent.

However if $\sigma(A)=\sigma_{s}(A)$ is more complicated there 
may exist two different rank one perturbations of $A$ which are unitarily
equivalent. 

To show that, it is enough to give an example of  two  equivalent pairs
of singular measures $(\mu,\nu)$ and $(\mu',\nu')$ 
with different phase shifts. Indeed, if such pairs exist then 
$A_\mu+(\cdot,1)1$ and
$A_\mu+(\cdot,\sqrt f)\sqrt f$, where 
 $f\in L^1(\mu), f\geq 0, \mu'=f\mu$, 
are unitarily equivalent and $f\neq const$.

To construct such an example one can notice, for instance, that function $u$ in the proof of Theorem 5.5
is not unique. Indeed, the reenumeration of the sequences $\Cal A$ and $\Cal B$ 
in part a)
of the proof of the implication i)$\Rightarrow$ ii) 
can be done in many different ways. It can be shown
that some of the enumerations will give us different phase shifts $u$ at the end of our
construction.  For instance, let  $n_k,m_k\in\Bbb N$ be such as
in part a)
of the proof of the implication i)$\Rightarrow$ ii). Then we can  choose $n'_k,m'_k\in\Bbb N$
to satisfy the conditions analogous to 1) and 2) and such that $n'_1=1, m'_k=1$ and
$|a_{n'_k}-b_{m'_k}|<|a_{n_k}-b_{m_k}|/2$ for $k=2,3,...$. Proceeding in the same way as 
in part a)
of the proof of the implication i)$\Rightarrow$ ii), we can define functions $u'_k$ and
consider the limit $u'$ of the sequence $\{u_k'\}$. Then 
$$
\left|\{\chi_{(a_1;b_1)}\neq u'\}\right|=\left|\{u'_1\neq u'\}\right|
<\left|\{u_1\neq u\}\right|/2=\left|\{\chi_{(a_1;b_1)}\neq u\}\right|.$$
Therefore $|\{u\neq u'\}|>0$. 
At the same time we again can prove that $u'$ is the phase shift of
the pair 
$(\sum \alpha_i'\delta_{a_i},\sum \beta_i'\delta_{b_i})$ 
for some positive $\alpha_i',\beta_i'$.
Therefore $u$ and $u'$ are different phase shifts of equivalent pairs of measures.

In terms of the families $M_\varphi$, we can say that if $(\mu,\nu)\in\Cal M_\varphi$
and $(\mu',\nu')\in\Cal M_\phi$ are equivalent pairs of measures and $\varphi$
is a finite Blacshke product, then $\varphi$ is equal to $b\circ\phi$ for some
M\"obius transform $b$. However, if $\varphi$ is more complicated then $\varphi$
and $\phi$ can be completely different functions.

In fact, modifying the construction from the proof of Theorem 5.5
one can show that there exist infinitely many different phase shifts
such that the corresponding pairs of measures are equivalent to
$(\sum \alpha_i\delta_{a_i},\sum \beta_i\delta_{b_i})$.
It means that
if $A$ and $B$ are two pure point operators whose eigenvalues are dense
in the same set $K\subset\Bbb R$ satisfying  condition ii) of Theorem 5.5,
then there exist infinitely many different (corresponding to noncollinear vectors) 
rank one perturbations of $A$
unitarily equivalent to $B$. 
However, it is absolutely unclear
if such situation is possible in the case when $A$ or $B$ has a nontrivial
singular continuous part.

\endremark

\heading 6. One example of the absence of the mixed spectrum \endheading

Let as usual $A$ be a cyclic self-adjoint operator, $\phi$  its
cyclic vector, $A_\lambda=A+\lambda(\cdot,\phi)\phi$ $\ \ \lambda\in
\Bbb R$. Let $\mu_\lambda$ be the spectral measure of  $\phi$
for $A_\lambda$. Let us denote by  $\Pi$ and $\Sigma$ the sets of $\lambda$ for which
$\mu_\lambda$ has  nontrivial pure point and nontrivial 
singular continuous part on $[0;1]$ respectively. Then the set $\Pi\cap\Sigma$
will consist of those $\lambda$ for which the corresponding measures are
``mixed'' on the interval $[0;1]$.

The  paper [D] of Donoghue gives examples in which 
$\Pi=\{0\},\Sigma=\Bbb R\setminus \{0\}$ and $\Sigma=\{0\},\Pi=\Bbb R\setminus \{0\}$.
Therefore the set $\Pi\cap\Sigma$ can be empty when both $\Pi$ and $\Sigma$ are nonempty.

One of the natural questions which arise from the 
recent results on rank one perturbations (see [R-J-L-S])
 is  {\it whether the set $\Pi\cap\Sigma$ can be empty (or almost empty) 
when the sets $\Pi$ and $\Sigma$
are sufficiently big (topologically or in measure)}.

The following example gives a partial answer to this question.

\example{Example 6.1}

We will show that there exist a self-adjoint cyclic operator $A$ such that
for some cyclic vector $\phi$ the operators $A+\lambda(\cdot,\phi)\phi$
are singular continuous on $[0;1]$ for all
$0\geq\lambda\geq 1$ 
and pure point for all other $\lambda\in\Bbb R$.

 We will first construct the  Krein spectral shift $u$.

Let $\{a_n\}$ be a sequence of real numbers  $0<a_n<1$ monotonically decreasing to 0
and such that
$$\prod_{n=1}^\infty (1-a_n)=c>0\tag6.1$$
and
$$1-\prod_{k=n}^\infty (1-a_k)\geq \frac1n.\tag6.2$$
Consider the Cantor set $C$ corresponding to the sequence $\{a_n\}$:
let  $$C_0=I_0^0=[0;1], \ \ C_1=I_1^1\cup I_2^1,...,\ \  C_n=I_1^n\cup...
\cup I_{2^n}^n,...$$ where $$I_{2k}^n\cup I_{2k-1}^n=I_k^{n-1}\setminus
\Delta_k^{n-1}$$
and $\Delta_k^n$ is the open interval placed in the center of the interval
$I_k^n$ and such that $|\Delta_k^n|=a_n|I_k^n|$ and let 
$C=\bigcap_{n=0}^\infty C_n$.

The Cantor set $C$ has the following properties:

 $$|C|=c\tag6.3$$
and
 $$\frac{|I_n^k\cap C|}{|I_n^k|}=
\prod_n^\infty(1-a_n)
\geq\frac1n\geq\frac{\ln 2}{-\ln |I_n^k|}\tag 6.4$$
for any $n,k\in\Bbb N$. 

Define $u=\pi$ on $C$ and $u=0$ elsewhere on $\Bbb R$.
Denote $U(z)=\Cal Ku(z)$ for $z\in\Bbb C_+$.

\proclaim{Claim} $U$ has a finite nontangential derivative $U'(x)=\lim_{
z \underset \ngtr \to \longrightarrow x} \frac {U(x)-U(z)}{x-z}$
at a point $x\in\Bbb R$ if and only if $x\not\in C$.
\endproclaim
\demo{Proof}
Since $u$ is locally constant on $\Bbb R\setminus C$, $U'$ obviously
exists there.

Let intervals $I^n_k$ be the same as in the construction of $C$ above.
Let $x\in C$. Let $\{I_{n_k}^k\}_{k=1}^\infty$ be the sequence of intervals
containing $x$. Denote by $x_k$ the middle of the interval $I_{n_k}^k$ and put
$y_k=|x-x_k|$. Then  (6.4)  imply that
$$\left|\pi -\Cal Pu(x_n+iy_n)\right|\to 0$$ (since $\prod_n^\infty(1-a_n)\to 0$)
but
$$\left|\pi -\Cal Pu(x_n+iy_n)\right|\geq\frac d{|\ln y_n|}$$
for some $d>0$.
Thus $U(x)-U(z)\neq O(x-z) $ as $z\underset\ngtr\to\rightarrow x$ and $U'(x)$ 
does not exist.$\blacktriangle$
\enddemo
 
Consider a real measure $\nu_0$ such that $u$ is the phase shift of a pair $(\nu_0,\nu)$
for some measure $\nu$.
Put $A=A_\nu, A_\lambda=A+\lambda(\cdot,1)1$. Let $\nu_\lambda$
be the spectral measure of $1$ for $A_\lambda$ (then $\nu=\nu_1$). 

By Lemma 3.2 $\nu$ is singular. Therefore all $\nu_\lambda$ are singular. Hence
$$|\Cal K\nu_\lambda|\nec\infty$$
for $\nu_\lambda$-a. e. $x$.
By formula (2.3) this means that
$$1+\pi\Cal K\nu_0=\exp U(z)\nec 1-\frac1{\lambda}$$ 
for $\nu_\lambda$-a. e. $x$.
Thus by the definition of the phase shift for all $\lambda\in (0;1)$
$$\arg(1+\pi\Cal K\nu_0(z))=\Cal Pu(z)\nec \pi$$
for $\nu_\lambda$-a. e. $x$. 
Thus all $\nu_\lambda, \lambda\in (0;1)$ are concentrated
on the set $C$. Since $U'$ does not exist on $C$
and $$\Cal K\nu_0=\exp U-1 \nec 1-1/\lambda\neq 0$$ $\nu_\lambda$-a. e., 
$(\Cal K\nu_0)'$ does not exist on $C$ for all $\lambda\in (0;1)$.
Hence the characteristic function $\phi$ defined in $\Bbb C_+$ by formula (2.4) does not have
a nontangential derivative on $C$. Thus
 by the Remark after Lemma 3.8
 all $\nu_\lambda, \lambda\in (0;\infty)$ are singular continuous.
Similarly all $\nu_\lambda$ for $ \lambda\in \Bbb R\setminus [0;1]$ are
concentrated on $\Bbb R\setminus C$. Since $U'$ exists everywhere on $\Bbb R\setminus C$, $\ \mu_\lambda$ for
$\lambda\in \Bbb R\setminus [0;1]$ are pure point.

To prove that $\nu_0$ is continuous, let us notice that since by Lemma 3.3
$$p.\ v.\ \int\limits_\Bbb R \frac{u(t+x)dt}{t}=\infty$$ for $\nu_0$-a.  e. $x$, $\nu_0$ is concentrated on $C$.

Let $y\in C$ and let $I^k_{n_k}$ be the sequence of intervals from the construction of $C$
containing $y$.
Then (6.4) implies that condition (3.8) is not satisfied at $y$.
Thus  $\nu_0$ can not have a point 
mass at $x$. Hence $\nu_0$ is continuous.
Similarly, $\nu_1$ is continuous.

\endexample
\remark{\bf Remark} It is still unclear if there exist examples such that
the sets $\Pi$ and $\Sigma$ are big (topologically or in measure) but $\Pi\cap\Sigma=\emptyset$ when
$\sigma(A)$ contains the interval $[0;1]$. Modifying the above example,
 we can obtain a dense in $[0;1]$ set $C'$ by inserting smaller Cantor sets
into each complimentary interval $I_n^k$, then inserting smaller Cantor sets 
 into each new complimentary interval
and so on. If the size of these Cantor sets  decreases to $0$
fast enough then replacing $C$ in the above example with $C'$ we will obtain
an example of $A$ and $\phi$ such that $\sigma(A)=[0;1]$, $A+\lambda(\cdot,\phi)\phi$
is continuous on $[0;1]$ {\it for  all} $\lambda\in [0;1]$ and pure point
{\it for almost all} $\lambda\in \Bbb R\setminus [0;1]$.

The case $\sigma(A)\supset [0;1]$ became especially interesting after it was shown in [Go] and [R-M-S] that
$\sigma(A)\supset [0;1]$ implies $\Sigma\setminus\Pi$ is a dense $G_\delta$.
\endremark


\Refs       
\refstyle{A}   
\widestnumber\key{R-J-L-S}
\ref\key A-C \by P.  Ahern and D. Clark \paper Radial limits and invariant
subspaces\jour Amer. J. Math.\vol 92\yr 1970\pages 332-342
\endref

\ref\key  A     \by  N. Aronszajn  \paper    
 On a problem of Weyl in the theory of singular Sturm-
Liouville equations
\jour Amer.  J.  Math.
\vol 79  \yr 1957  \pages 597-610
\endref

\ref\key A-D \by N. Aronszajn and W. Donoghue \paper
On exponential representations of analytic functions in the upper
half-plane with positive imaginary part \jour J. Analyse Math. \vol 5
 \yr 1956-1957 \pages 321-388 \endref
                 


\ref\key A1     
\by  A. B. Aleksandrov  \paper    
Multiplicity of boundary values of inner functions
\jour 
Izv.  Acad.  Nauk.  Arm.  SSR, Matematica 22\vol 5  \yr 1987 
\lang Russian
\pages 490-503
\endref

\ref\key  A2     
\by  A. B.  Aleksandrov  \paper    
Inner functions and related spaces of pseudo continuable 
functions 
\jour 
Proceedings of LOMI seminars 
\vol 170  \yr 1989 \lang Russian \pages 7-33
\endref

\ref\key  A3  \by A. B. Aleksandrov \book Private communications \endref


\ref\key C-P \by R. W. Carey and J. D. Pincus \paper
Unitary equivalence modulo the trace class for self-adjoint operators
\jour Amer J. Math. \vol 98 \yr 1976 \pages 481-514 \endref

\ref\key  C     \by D. Clark  \paper    
One dimensional perturbations of restricted shifts
\jour 
J.  anal.  math.  
\vol 25  \yr 1972  \pages 169-91
\endref


\ref\key  D     \by W. Donoghue  \paper    
On the perturbation of spectra 
\jour 
 Comm.  Pure Appl.  Math. 
\vol 18  \yr 1965  \pages  559-576
\endref
\ref \key G  \by J. B. Garnett \book Bounded analytic functions
\bookinfo Academic Press \publaddr New York \yr 1981 \endref

\ref \key Go \by A. Gordon \paper Pure point spectrum under 1-parameter
perturbations and instability of Anderson localization
\jour to appear in Commun. Math. Phys. \endref

\ref\key K \by M. G. Krein \paper Perturbation determinants and a
formula for the traces of unitary and self-adjoint operators
\jour Sov. Math. Dokl. \vol 3 \yr 1962 \pages 707-710 \endref

\ref\key M-P \by M. Martin and M. Putinar \book Lectures on Hyponormal
operators \publ Operator Theory: Advances and Applications \vol 39 \yr 1989 \endref

\ref\key P1 \by A. Poltoratski \paper On the boundary behavior of pseudocontinuable
functions \jour St. Petersburg Math. J.\vol  5 \yr 1994 \pages 389-406\endref

\ref\key P2 \by A. Poltoratski \paper On the distributions of boundary values of Cauchy integrals
\jour to appear in Proc. Amer. Math. Soc.\endref

\ref\key  R-J-L-S  \by R. del Rio, S. Jitomirskaya, Y. Last
and B. Simon 
\paper Operators with singular continuous spectrum, 4. 
Hausdorff dimension and rank one perturbations
\jour to appear
\endref

\ref \key R-M-S \by R. del Rio, N. Makarov and B. Simon \paper
Operators with singular continuous spectrum, II. Rank one operators
\jour to appear in Commun. Math. Phys. \endref 

\ref\key  S     \by B. Simon  \paper    
Spectral analysis of rank one perturbations and 
applications
\jour 
Proc. 1993 Vancouver Summer School in Mathematical Physics
\endref

\ref\key  Sa  \by D.~Sarason \book Sub-Hardy Hilbert Spaces in the unit disk
 \bookinfo The University of Arkansas lecture notes in the mathematical sciences; v. 10
\publaddr New York \yr 1994 \publ J. Wiley and  Sons
\endref


             
\ref\key   S-W    
\by B. Simon, T. Wolff  \paper    
 Singular continuous spectrum under
rank one perturbations and localization for random Hamiltonians
\jour 
Comm.  Pure
Appl.  Math.
\vol  39 \yr 1986  \pages 75-90
\endref



\endRefs
\enddocument